\documentclass[letterpaper, 10 pt, conference]{ieeeconf}  

\IEEEoverridecommandlockouts                              

\usepackage{amsmath,amssymb,amsfonts}

\UseRawInputEncoding
\usepackage{caption}
\usepackage{subcaption}
\def\disp{\displaystyle}
\def\e{\epsilon}
\def\vt{\vartheta}

\def\gph{\mbox{\rm gph}\,}
\def\epi{\mbox{\rm epi}\,}
\def\dom{\mbox{\rm dom}\,}
\def\dn{\downarrow}

\def\ox{\bar{x}}
\def\oy{\bar{y}}

\def\ou{\bar{u}}

\def\oe{\bar{e}}
\def\({\left(}
\def\){\right)}
\def\[{\left[}
\def\]{\right]}

\def\gg{\gamma}

\def\gph{\mbox{\rm gph}\,}
\def\epi{\mbox{\rm epi}\,}

\def\n{\Big\Vert}
\def\en{\Big\Vert}

\def\la{\langle}
\def\ra{\rangle}
\def\vph{\varphi}
\def\emp{\emptyset}
\def\ph{\varphi}
\def\be{\beta}
\def\lm{\lambda}

\def\tto{\rightrightarrows}
\def\th{\theta}
\def\O{\Omega}

\def\oR{\Bar{\R}}

\newcommand{\R}{\mathbb{R}}
\newcommand{\N}{\mathbb{N}}

\newtheorem{theorem}{Theorem}[section]

\newtheorem{definition}[theorem]{Definition}

\usepackage{graphicx,color}

\include{Xlatin1}

\newcommand{\dist}{\mathop{\rm dist}}

\newcommand{\ba}{\begin{array}}
\newcommand{\ea}{\end{array}}

\title{\LARGE \bf Optimization of Controlled Free-Time Sweeping Processes with Applications to Marine Surface Vehicle Modeling}

\author{Tan H. Cao$^{1}$, Nathalie T. Khalil$^{2}$, Boris S. Mordukhovich$^{3}$,\\ Dao Nguyen$^{3}$, Trang Nguyen$^{3}$, Fernando Lobo Pereira$^{2}$
\thanks{$^{1}$T.H. Cao is with Department of Applied Mathematics and Statistics, SUNY (State University of New York) Korea
    {\tt\small tan.cao@stonybrook.edu}}%
\thanks{$^{2}$N.T. Khalil and F.L. Pereira are with SYSTEC, Faculty of Electrical Engineering, Porto University, and with the Institute for Systems and Robotics, 4200-465 Porto, Portugal
        {\tt\small khalil.t.nathalie@gmail.com, flp@fe.up.pt}}%
    \thanks{$^{3}$B.S. Mordukhovich, D. Nguyen, and T. Nguyen are with Department of Mathematics, Wayne State University, USA
        {\tt\small boris@math.wayne.edu, dao.nguyen2@wayne.edu, daitrang.nguyen@wayne.edu}}%
}
\begin{document}

\maketitle
\thispagestyle{empty}
\pagestyle{empty}

\begin{abstract}
The paper is devoted to a free-time optimal control problem for sweeping processes. We develop a constructive finite-difference approximation procedure that allows us to establish necessary optimality conditions for discrete optimal solutions and then show how these optimality conditions are applied to solving a controlled marine surface vehicle model.
\end{abstract}

\section{Introduction and Problem Formulation}\label{sec:intro}

The {\em sweeping process} was introduced by Moreau in the 1970s (see \cite{moreau}) in the form of the differential inclusion form
\begin{equation*}
\left\{\begin{matrix}
\dot{x}(t)\in-N\big(x(t);C\big)\;\textrm{ a.e. }\;t\in[0,T],\\x(0)=x_0\in C\subset\R^n.
\end{matrix}\right.
\end{equation*}
The sweeping process and its modifications have been extensively studied and applied to various fields including aerospace engineering, process control, robotics, bioengineering, chemistry, biology, economics, finance, management science, and engineering. Moreover, the sweeping dynamics play an important role in the theory of variational inequalities and complementarity problems. Among applications in mechanical and electrical engineering, we mention mechanical impact, Coulomb friction, diodes and transistors, queues and resource limits, etc.; see, e.g., the recent survey in \cite{bt}.

Optimal control problems for various types of sweeping processes have been formulated much more recently (see \cite{cmn18a} and the references therein), while being realized as very challenging control theory due to high discontinuity of the controlled sweeping dynamics and the unavoidable presence of hard state constraints. Nevertheless, within a rather short period of time, many important results have been obtained on necessary optimality conditions for controlled sweeping processes with valuable applications to friction and plasticity, robotics, traffic equilibria, ferromagnetism, hysteresis, economics, and other fields of engineering and applied sciences; see, e.g., \cite{ac,b1,cmn18a,pfs,zeidan} with more references and discussions. Let us mention to this end the recent papers \cite{vbp1,vbp2}, where optimal control problems for {\em linear complementarity systems} have been studied and applied to practical models that are highly important in the area of {\em Automatic Control}. Such problems can be written in a form of controlled sweeping processes, where $C$ is an orthant in $\R^n$. However, there are many great unsolved problems in optimal control theory for sweeping processes with strong requirements for further applications. Some of these issues, from both viewpoints of theory and applications, are addressed in this paper.

Here, we consider the following free-time optimal control problem labeled as $(P)$:
minimize the Mayer-type cost functional which depends explicitly on the final time
\begin{equation*}
J[x,u,T]:=\varphi(x(T),T)
\end{equation*}
over control functions $u(\cdot)$ and the corresponding trajectories $x(\cdot)$ satisfying the system
\begin{equation}\label{Problem}
\left\{\begin{matrix}
\dot{x}(t)\in-N\big(x(t);C\big)+g\big(x(t),u(t)\big)\;\textrm{ a.e. }\;t\in[0,T],\\x(0)=x_0\in C\subset\R^n,\\
u(t)\in U\subset\R^d\;\textrm{ a.e. }\;t\in[0,T],
\end{matrix}\right.
\end{equation}
where the set $C$ is a convex polyhedron given by
\begin{equation}\label{C}
\left\{\begin{matrix}
C:=\bigcap_{j=1}^{s}C^j\textrm{ with }C^j:=\{ x\in\R^n\big|\;\la x^j_*,x\ra\le c_j\},\\
\|x_\ast^j\|=1,\;j=1,\ldots,s,
\end{matrix}\right.
\end{equation}
with $N(x;C)$ standing for the normal cone of convex analysis. From \eqref{Problem} we automatically have the {\em state constraints}
\begin{eqnarray*}\label{e:8}
x(t)\in C,\textrm{ i.e., }\la x^j_*,x(t)\ra\le c_j\;\textrm{ for all }\;t\in [0,T]\;\nonumber\\\textrm{ (with different } T)\;\textrm{ and }\;j=1,\ldots,s.
\end{eqnarray*}
Note that defining state constraints implicitly via the domain of the normal cone in (1) is significantly different from the formulation of pure state constraints in standard control theory. In what follows, we identify the arc $x:[0,T]\to\R^n$ with its extension to $(0,\infty)$ defined by
\begin{equation*}
x_e(t):=x(T) \textrm{ for all }t>T.
\end{equation*}
Given $x(\cdot)\in W^{1,2}([0,T],\R^n)$ with the norm
\begin{equation*}
\n x\en_{W^{1,2}}:=\n x(0)\en+\n \dot{x}_e\en_{L_2},
\end{equation*}
we specify the notion of local minimizers studied below. For simplicity, suppose that the set $g(x;U)$ is convex, which actually does not much restrict the generality; cf. \cite{cmn18a}.
\begin{definition}\label{Def3.1}
A feasible solution $(\ox(\cdot),\ou(\cdot),\overline{T})$ for $(P)$ is a {\em $W^{1,2}\times L^2$-local minimizer} to this problem if there exists $\e>0$ such that $J[\ox,\ou,\overline{T}]\le J[x,u,T]$ for all feasible solutions $(x(\cdot),u(\cdot),T)$ satisfying the constrains of $(P)$ and
\begin{equation*}
\int_0^{\overline{T}}\(\n\dot{x}_e(t)-\dot{\ox}_e(t)\en^2+\n u(t)-\ou(t)\en^2\)dt+(\overline{T}-T)^2<\e.
\end{equation*}
\end{definition}\vspace*{0.05in}

Our approach to investigate problem $(P)$ in order to establish necessary optimality conditions for its local minimizers is based on the {\em method of discrete approximations} developed in \cite{m95,Mord} for Lipschitzian differential inclusions and then extended in \cite{b1,CaoMordukhovich2018DCDS,chhm,cmn18b,cmn18a,pfs} to various kinds of controlled sweeping processes. This method consists of constructing well-posed discrete approximations of $(P)$ whose solutions strongly converge to the prescribed local minimizers of $(P)$, then deriving necessary optimality conditions for discrete-time problems, and finally establishing by passing to the limit, with the discretization step decreasingly converging to zero, necessary optimality conditions for local minimizers of $(P)$. Due to the text size limitation, we present here only optimality conditions for discrete approximations, which give us sufficient information to solve some applied optimal control problems that arising in marine surface vehicle modeling and control.

The rest of the paper is organized as follows. In Section~\ref{sec:assumptions} we present the standing assumptions on the problem data. Section~\ref{sec:va} reviews the tools of variational analysis used below. Section~\ref{sec:NCO} is devoted to the formulation of  necessary optimality conditions for discrete approximations of $(P)$. A proof outline with the key ideas is given in Section~\ref{sec:proof}. Section~\ref{sec:examples} contains applications of the obtained necessary optimality condition to the controlled marine surface vehicle model with providing numerical calculations in typical settings for two marine surface vehicles. We finish with Section~\ref{sec:Conclusions} that contains concluding remarks and discussions of some topics of our future research.

\section{Standing Assumptions}\label{sec:assumptions}
First we present the following standing assumptions:
\begin{itemize}
\item[(H1)] The set $U\ne\emp$ is closed and bounded in $\R^d$.
\item[(H2)] The {\em positive linear independence constraint qualification} $($PLICQ$)$ holds at $x(t)$ on $[0,T]$ with varying $T$:
\begin{eqnarray*}\label{plicq}
\Big[\sum_{j\in I(x)}\lambda_jx^j_\ast=0,\;\lambda_j\in\R_+\Big]\nonumber\\\Longrightarrow\big[\lambda_j=0\textrm{ for all }j\in I(x)\big],
\end{eqnarray*}
where the active index set $I(x)$, $x\in C$, is defined by
\begin{equation}\label{aci}
I(x):=\big\{j\in\{1,\ldots,s\}\;\big|\;\la x^j_\ast,x\ra= c_j\big\}.
\end{equation}
\item[(H3)] The perturbation mapping $g\colon\R^n\times U\to\R^n$ is  Lipschitz continuous with respect to $x$ uniformly on $U$ whenever $x$ belongs to a bounded subset of $\R^n$ and satisfies there the sublinear growth condition
\begin{equation*}
\|g(x,u)\|\le\be\big(1+\|x\|\big)\;\mbox{ for all }\;u\in U
\end{equation*}
with some positive constant $\be$.
\end{itemize}\vspace*{0.03in}
Define the set-valued mapping $F\colon\R^n\times\R^d\tto\R^n$ by
\begin{equation}\label{F0}
F(x,u):=N(x;C)-g(x,u)
\end{equation}
and deduce from the theorem of the alternative that
\begin{equation}\label{F}
F(x,u)= \Big\{\sum_{i\in I(x)}\lm^i x^i_*\;\Big|\;\lm^i\ge 0\Big\}-g(x,u).
\end{equation}

\section{Tools of Variational Analysis}\label{sec:va}

Let us recall the tools of variational analysis employed below; see \cite{Mord} and \cite{rw} for more details. The (Painlev\'e-Kuratowski) {\em outer limit} of a set-valued
mapping/multifunction $F\colon\R^n \tto\R^m$ at $\ox$ with $F(\ox)\ne\emp$ is 
\begin{align*}\label{Pa}
\underset{x\to\ox}{\textrm{Lim sup }}F(x):=\big\{y\in\R^m\;\big|\;\exists\textrm{ sequences }\;x_k\to\ox,\\
y_k\to y\textrm{ such that }\;y_k\in F(x_k),\;k\in\N\big\}.
\end{align*}
The (basic, limiting, Mordukhovich) {\em normal cone} to subsets $\O\subset\R^n$ that are locally closed around $\ox$ is given by
\begin{equation}\label{nor}
N(\ox;\O)=N_\O(\ox):=\underset{x\to\ox}{\textrm{Lim sup}}\big\{\textrm{cone}[x-\Pi(x;\O)]\big\},
\end{equation}
where $\Pi(x;\O):=\big\{u\in\O\;\big|\;\|x-u\|=\dist(x;\O)\big\}$ is the Euclidean projection of $x$ to $\O$, and where ``cone" stands for the (generally nonconvex) conic hull of the set.

We consider a set-valued mapping $F\colon\R^n\tto\R^m$ locally closed around a graph point $(\ox,\oy)\in\gph F$. The {\em coderivative} of $F$ at $(\ox,\oy)$ is defined by
\begin{equation*}
D^*F(\ox,\oy)(u):=\big\{v\in\R^n\;\big|\;(v,-u)\in N\big((\ox,\oy);\gph F\big)\big\},
\end{equation*}
where $u\in\R^m$. If $F\colon\R^n\to\R^m$ is single-valued and smooth around $\ox$, then we have
\begin{equation*}
D^*F(\ox)(u)=\big\{\nabla F(\ox)^*u \big\}\;\textrm{ for all }\;u\in\R^m,
\end{equation*}
where $\nabla F(\ox)^*$ is the adjoint/transposed Jacobian matrix, and where $\oy=F(\ox)$ is omitted.

Given an extended-real-valued and lower semicontinuous function $\varphi\colon\R^n\to\oR:=(-\infty,\infty]$ with $\ph(\ox)<\infty$, the {\em subdifferential} of $\ph$ at $\ox\in\dom\varphi$ is defined via the normal cone \eqref{nor} to the epigraph of $\ph$ by
\begin{equation*}
\partial\varphi(\ox):=\big\{v\in\R^m\;\big|\;(v,-1)\in N\big((\ox,\varphi(\ox));\epi\varphi\big)\big\}.
\end{equation*}
These normal cone, coderivative, and subdifferential enjoy {\em full calculus}, which can be found in \cite{Mord}, and \cite{rw}.

\section{Necessary Optimality Conditions}\label{sec:NCO}
We first construct a sequence of discrete approximation problems with a varying grid whose optimal solutions {\em strongly} converge in $W^{1,2}$-norm to a prescribed local minimizer $\{\bar{x}(\cdot),\bar{u}(\cdot),\overline{T}\}$ of the original problem $(P)$. For simplicity we replace the derivative $\dot x(t)$ in \eqref{Problem} by
\begin{equation*}
\dot x(t)\approx\frac{x(t+h)-x(t)}{h}\textrm{ as }h\dn 0.
\end{equation*}
Whenever $k\in\N$, take $T_k$ close to $T$ and set the grid
\begin{equation}\label{grid}
\left\{\begin{matrix}
t^k_0=0,\;t^k_k=T_k,\\
t^k_{i+1}=t^k_i+h^k_i,\;i=0,\ldots,k-1.
\end{matrix}\right.
\end{equation}
Let $(\ox(\cdot),\ou(\cdot),\overline{T})$ be a $W^{1,2}\times L^2$-local minimizer of problem $(P)$ for the differential inclusion \eqref{Problem}. Define the {\em discrete approximation problem} $(P_{k})$ by:
\begin{align*}
&\mbox{minimize}\;\;J_k[x^k,u^k,T_k]:=\varphi(x^k_k,T_k)+(T_k-\overline T)^2+\nonumber\\
&\sum_{i=0}^{k-1}\int_{t^k_i}^{t^k_{i+1}}
\left(\n\frac{x^k_{i+1}-x^k_i}{h^k_i}-\dot{\bar{x}}(t)\en^2+ \right.\nonumber
\left. \n u^k_i-\bar{u}(t)\en^2\right) dt
\end{align*}
over $(x^k,u^k,T_k):=(x^k_0,x^k_1,\ldots,x^k_{k-1},u^k_0,u^k_1,\ldots,u^k_{k-1},T_k)$ satisfying the following constraints:
\begin{equation*}
x^k_{i+1}- x^k_i\in -h^k_i F(x^k_i,u^k_i)\;
\textrm{ for }\;i=0,\ldots,k-1,
\end{equation*}
\begin{equation*}
x^k_0:=\ox_0\in C,
\;u^k_0:=\ou (0),
\end{equation*}
\begin{equation*}
\sum_{i=0}^{k-1}\int_{t^k_i}^{t^k_{i+1}}\(\n\frac{x^k_{i+1}-x^k_i}{h^k_i}-\dot{\ox}(t)\en^2+\n u^k_i-\ou(t)\en^2\)dt\le \e,
\end{equation*}
\begin{equation*}
u^k_i\in U\;\textrm{ for }\;i=0,\ldots,k-1,\;\;|T_k-\overline{T}|\le\e,
\end{equation*}
\begin{equation*}
\n\(x^k_i,u^k_i\)-\(\ox(t^k_i),\ou(t^k_i)\)\en\le\e\;\textrm{ for }i=0,\ldots,k-1,
\end{equation*}
\begin{equation}\label{re_7}
\la x^j_\ast, x^k_k\ra \leq c_j \textrm{ for all } j =1,\ldots,s,
\end{equation}
where $\e >0$ as in Definition~\ref{Def3.1}. In each problem $(P_k)$, the final time $T_k$ and the discretization step $h^k_i$ are \textit{variable} for any fixed $k\in\N$.\vspace*{0.03in}

Now, we present necessary optimality conditions for $(P_k)$ expressed entirely via the given data. Having $I(x)$ from \eqref{aci} and $y\in\R^n$, consider $j\in I(x)$ for which
\begin{align*}
I_0(y):=\{j|\la x^j_*,y\ra=c_j\}\;\mbox{ and }\;I_>(y):=\{j|\la x^j_*,y\ra>c_j\}.
\end{align*}

\begin{theorem}\label{theorem:main thm} Let $(\ox^k(\cdot),\ou^k(\cdot),\overline{T}_k)$ be an optimal solution to problem $(P_k)$, where the cost function $\ph$ is locally Lipschitzian around $(\ox^k(\overline{T}_k), \overline{T}_k)$ in addition to the standing assumptions. Then, there exist dual elements $(\mu^k_0,q^k,p^k)$ together with vectors $\eta^k_i\in\R^s_+$ for $i=0,\ldots,k$, and $\gg^k_i\in\R^s$ for $i=0,\ldots,k-1$ satisfying the following:
\begin{itemize}
\item[1.] The {\sc nontriviality conditions}. In general we have
\begin{equation*}
\mu^k_0+\n\eta^k_k\en+\sum_{i=0}^{k-1}\n p^k_i\en+\n q^k\en\ne 0.
\end{equation*}
If the matrices $\nabla_u g(\ox^k_i,\ou^k_i)$ are of full rank as $i=0,\ldots,k-1$, the {\em enhanced nontriviality condition} holds:
\begin{equation*}
\mu^k_0+\|\eta^k_k\|+\|p^k_0\|+\|q^k\|\ne 0.
\end{equation*}
\item[2.] The {\sc primal-dual dynamic relationships}:\\
$\bullet$ The {\em primal arc representation}
\begin{equation*}
-\frac{\ox^k_{i+1}-\ox^k_i}{h^k_i}+g(\ox^k_i,\ou^k_i)=\sum_{j\in I(\ox^k_i)}\eta^k_{ij} x^j_*.
\end{equation*}
$\bullet$ The {\em adjoint dynamic systems}
\begin{align*}
&\hspace{-.5cm}\dfrac{p^k_{i+1}-p^k_i}{h^k_i}=-\nabla_x g(\ox^k_i,\ou^k_i)^*\Big(-\dfrac{\mu^k_0\xi^k_{iy}}{h^k_i}+p^k_{i+1}\Big)+\nonumber\\
&\hspace{2.5cm}\disp\sum_{j\in\mathcal I }\gg^k_{ij} x^j_*,\nonumber\\
&\hspace{-.5cm}\textrm{where  }\nonumber\\
&\hspace{-.5cm}\mathcal I:= I_0\(-\frac{\mu^k_0\xi^k_{iy}}{h^k_i}+p^k_{i+1}\)\cup I_>\(-\frac{\mu^k_0\xi^k_{iy}}{h^k_i}+p^k_{i+1}\).
\end{align*}
$\bullet$ The {\em adjoint inclusions}
\begin{align*}
&\hspace{-.5cm}\(-p^{k}_k-\disp\sum^{s}_{j=1}\eta^k_{kj}x^j_*,\;\bar
{H}^{k}+2\mu^k_0(\overline{T}-\overline{T}_{k})+\mu^k_0\varrho_k\)\in\nonumber\\
&\hspace{1.5cm}\partial\Big(\mu^k_0\varphi\Big)(\bar{x}^{k}(\overline{T}_{k}),\overline{T}_{k}),\nonumber
\end{align*}
where we use the notation
\begin{align*}
&\bar{H}^k:=\dfrac{1}{k}\sum_{i=0}^{k-1}\la p^k_{i+1},y^k_i\ra,\\
&\varrho_k:=\sum_{i=0}^{k-1}\left[\dfrac{i}{k}\n\(\frac{\ox^k_{i+1}-\ox^k_i}{h^k_i}-\dot{\ox}(t_i),\ou^k_i-\ou(t_i)\)\en^2\right.\\
&\left. -\dfrac{i+1}{k}\n\(\frac{\ox^k_{i+1}-\ox^k_i}{h^k_i}-\dot{\ox}(t_{i+1}),\ou^k_i-\ou(t_{i+1})\)\en^2\right],\nonumber\\
&\xi^k_i \nonumber =\(\xi^k_{iu},\xi^k_{iy}\)\textrm{ with }\xi^k_{iu}:=\int_{t^k_i}^{t^k_{i+1}}
\(\ou^k_i-\ou(t)\)dt\\
&\textrm{ and }\;\xi^k_{iy}:=\int_{t^k_i}^{t^k_{i+1}}\(\frac{\ox^k_{i+1}-\ox^k_i}{h^k_i}-\dot{\ox}(t)\)dt.
\end{align*}
\item[3.] The {\sc local maximum principle:}
$q^k_i\in N(\ou^k_i;U)$ as $i=0,\ldots,k-1$ with
\begin{equation*}
-\dfrac{\mu^k_0\xi^k_{iu}}{h^k_i}-\dfrac{q^k_i}{h^k_i}=-\nabla_u g(\ox^k_i,\ou^k_i)^*\Big(-\dfrac{\mu^k_0\xi^k_{iy}}{h^k_i}+p^k_{i+1}\Big),
\end{equation*}
which yields the {\em linearized global} form of the maximum principle when $U$ is convex.
\item[4.] The {\sc complementarity slackness conditions.} The following implications hold:
\begin{equation*}
\Big[\la x^j_\ast,\ox^k_i\ra<c_j\Big]\Longrightarrow\eta_{ij}^k=0,
\end{equation*}
\vspace{-.5cm}
\begin{eqnarray*}
\left\{\begin{matrix}
\Big[j\in I_>(-\frac{\mu^k_0\xi^k_{iy}}{h^k_i}+p^k_{i+1})\Big]\Longrightarrow\gg^k_{ij}\ge 0,\nonumber\\
\Big[j\notin \mathcal I\Big]\Longrightarrow\gg^k_{ij}=0,
\end{matrix}\right.
\end{eqnarray*}
\begin{equation*}
\[\la x^j_*,\ox^k_i\ra<c_j\]\Longrightarrow\gg^k_{ij}=0,
\end{equation*}
\begin{equation*}
\big[\la x^j_\ast,\ox^k_k\ra < c_j \big]\Longrightarrow\eta^k_{kj}=0,
\end{equation*}
where $i=0,\ldots,k-1$ and $j=1,\ldots,s$.\\
Finally, imposing the linear independence of the vectors $\{x^j_*|\;j\in I(\ox^k_i)\}$ ensures the implication
\begin{equation*}
\eta^k_{ij}>0\Longrightarrow\Big[\Big\la x^j_*,-\frac{\mu^k_0\xi^k_{iy}}{h^k_i}+p^k_{i+1}\Big\ra=0\Big].
\end{equation*}
\end{itemize}
\end{theorem}

\section{Brief Outline of Key Ideas of the Proof}\label{sec:proof}
For any fixed $k\in \N$ and $\e>0$, consider the following problem $(MP)$ with respect to variables $z:=(x^k_0,\ldots,x^k_k,u^k_0,\ldots,u^k_{k-1},y^k_0,\ldots,y^k_{k-1},\th)$:
\begin{align*}
\textrm{minimize }\;\phi_0(z):=\vph\big(x^k_k,\th\big)+(\th-\overline{T})^2+\\
\sum_{i=0}^{k-1}\int_{\frac{i\th}{k}}^{\frac{(i+1)\th}{k}}\n\(y^k_i-\dot{\ox}(t),u^k_i-\ou(t)\)\en^2dt\;\mbox{ s.t. }
\end{align*}
\begin{align*}
&\phi(z):=\sum_{i=0}^{k-1}\int_{\frac{i\th}{k}}^{\frac{(i+1)\th}{k}}\n\(y^k_i-
\dot{\ox}(t),u^k_i-\ou(t)\)\en^2dt -\e\le 0,\\
&\phi_j(z):=\la x^{j}_\ast,x^k_k\ra-c_{j}\le 0,\\
&g_i(z):=x^k_{i+1}-x^k_i-\dfrac{\th}{k} y^k_i=0,\\
&\Xi_i:=\big\{(x^k_0,\ldots,y^k_{k-1},\th)\;\big|\;-y^k_i\in F\(x^k_i,u^k_i\)\big\},\\
&z\in\Xi_k:=\big\{(x^k_0,\ldots,y^k_{k-1},\th)\;\big|\;x^k_0\;\textrm{ is fixed},\;\th \in \mathbb{R}_{+}\},\\
&\Xi'_i:=\left\{(x^k_0,\ldots,y^k_{k-1},\th)\;\Big|\; u^k_i\in U\right\}
\end{align*}
for all $i=0,\ldots,k-1,\;j=1,\ldots,s$. The necessary optimality conditions for $(MP)$ are given in \cite[Theorem~5.21(iii)]{Mord}. It follows from them that for  all $i=0,\ldots,k-1$ we have
\begin{align*}
&\Big(\frac{p^k_{i+1}-p^k_i}{h^k_i},-\frac{\mu^k_0\xi^k_{iu}}{h^k_i}-\frac{q^k_i}{h^k_i}\Big)\in\nonumber\\
&D^*F\Big(\ox^k_i,\ou^k_i,-\frac{\ox^k_{i+1}-\ox^k_i}{h^k_i}\Big)\(-\frac{\mu^k_0\xi^k_{iy}}{h^k_i}+p^k_{i+1}\).
\end{align*}
Taking into account the structures of $F$ in \eqref{F0}, \eqref{F} and using the coderivative calculations from \cite[Theorem~5.1]{cmn18a}, we arrive at all the conditions of
Theorem~\ref{theorem:main thm}.

Passing to the limit as $k\to\infty$ leads us to necessary optimality conditions for the original problem $(P)$.

\section{Marine Surface Vehicle Model}\label{sec:examples}
In this section we show how the obtained necessary optimality conditions allow us to find optimal solutions for a controlled marine surface vehicle model. The model deals with $n\ge 2$ unmanned surface vehicles (USVs) that have arbitrary shapes navigating at sea identified as virtual disks of different radii $R_i$, $i=1,\ldots,n$ on the plane. Each USV aims at reaching the target by the shortest pass with the minimum time $T$ while avoiding the other $n-1$ static and/or dynamic obstacles. The dynamics of this model are described in \cite{hb1} as an uncontrolled sweeping process.

Define the configuration space of USVs at time $t$ by $x=(x^1,\ldots,x^n)\in\R^{2n}$, where $x^i=(\|x^i\|\cos\th_i^x,\|x^i\|\sin\th_i^x)\in\R^2$ denotes the Cartesian position of the $i$-th vehicle, $\th_i^x$ stands for the constant direction that is the smallest positive angle in standard position formed by the positive $x$-axis and vectors $Ox^i$ with $0\in\R^2$ as the target.

The configuration is admissible when the motion of different USVs is safe by imposing the noncollision or nonoverlapping condition. This can be formulated mathematically as
\begin{eqnarray*}
A:=\big\{x=\(x^1,\ldots,x^n\)\in\mathbb{R}^{2n}\big|\;D_{ij}(x)\ge 0\big\}
\end{eqnarray*}
for all $i,j\in\{1,\ldots,n\}$, where $D_{ij}(x):=\|x^{i}-x^j\|-(R_i+R_j)$ is the distance between the disks $i$ and $j$.

The safe navigation of USVs can be described as follows. Starting from an admissible configuration at time $t_k\in[0,T]$ (with different $T$), consider $x_k:=x(t_k)\in A$, and then, get the next configuration after the period of time $h>0$, as $x_{k+1}=x(t_k+h)$. To ensure a safe navigation of all USVs at $t_k+h$ for a small value of $h>0$, the next configuration should also be admissible, i.e., $x(t_k+h)\in A$. This implies that the constraint $D_{ij}\(x(t_k+h)\)\ge 0$ should be satisfied. To verify this, employ the first-order Taylor expansion at $x_k\ne 0$ and reduce the constraint on the velocity vector to
\begin{eqnarray*}
D_{ij}\big(x(t_k+h)\big)=D_{ij}(x(t_k))+h\nabla D_{ij}(x(t_k))\dot{x}(t_k)+o(h)
\end{eqnarray*}
for small $h>0$. This constraint will be used to construct the next configuration in order to avoid the collision with static and/or dynamic obstacles. In this regard, set the vector $V(x)$ to be the desired velocity of all USVs. The admissible velocities preventing collisions during the navigation of USVs are defined as
\begin{eqnarray*}
&C_h(x):=\big\{V(x)\in\R^{2n}\big|\;D_{ij}(x)+h\nabla D_{ij}(x)V(x)\ge 0\big\}
\end{eqnarray*}
for all $i,j\in\{1,\ldots,n\},\;i<j$ and $x\in\R^{2n}$. Taking the admissible velocity $\dot{x}(t_k)\in C_h(x_k)$ gives us
\begin{eqnarray*}
D_{ij}(x_k)+h\big\la\nabla D_{ij}(x_k),\dot{x}(t_k)\big\ra\ge 0.
\end{eqnarray*}
Skipping the term $o(h)$ for small $h$, we deduce that $D_{ij}(x(t_k+h))\ge 0$, i.e., $x(t_k+h)\in A$.

In the absence of obstacles, the desired velocities are given by $V(x)=g(x)\in\R^{2n}$. In the presence of obstacles, the algorithm in \cite{hb1} seeks for optimal velocities to escape from surrounding obstacles by solving the following convex constrained optimization problem:
\begin{eqnarray}\label{Pi}
\mbox{minimize }\;\|g(x,u)-V(x)\|^2\;\mbox{ s.t. }\;V(x)\in C_{h}(x),
\end{eqnarray}
where the control $u$ is involved into the desired velocity term to adjust the actual velocities of the USVs and make sure that they do not overlap. The velocities can be modeled as
\begin{eqnarray*}\label{t:99**}
&g\big(x(t),u(t)\big)=\big(s_1\|u^1(t)\|\cos\th_1^u(t),s_1\|u^1(t)\|\sin\th_1^u(t),\nonumber\\
&\ldots,s_n\|u^n(t)\|\cos\th_n^u(t),s_n\|u^n(t)\|\sin\th_n^u(t)\big),
\end{eqnarray*}
where $s_i$ denotes the speed of the USV $i$, and $\th^u_i$ stands for the corresponding constant direction which is the smallest positive angle in standard position formed by the positive $x$-axis and vectors $u^i(t)$, with practically motivated control constraints represented by
\begin{eqnarray}\label{t:99*}
u(t)=\big(u^1(t),\ldots,u^n(t)\big)\in U\;\mbox{ for a.e. }t\in[0,T],
\end{eqnarray}
where the control set $U\subset\R^n$ will be specified below.

The algorithmic design in \eqref{Pi} means that $V_{k+1}$ is selected as the unique element from the set of admissible velocities as the one closest to the desired velocity $g(x,u)$ while avoiding overlapping. Consequently, the proposed scheme seeks for new directions $V(x)$ of USVs close to the desired direction $g(x,u)$ in order to bypass the surrounding obstacles. The desired position of the next configuration of marine vessels is generated as $x_{\textrm{ref}}(t+h)=(x_{\textrm{ref}},y_{\textrm{ref}})=x(t)+hV(x)$ and the desired via-point posture position of the marine craft is
$$
\eta_{\textrm{ref}}=(x_{\textrm{ref}},y_{\textrm{ref}},\psi_{\textrm{ref}}),\quad\psi_{\textrm{ref}}=\tan^{-1}\(\dfrac{y_{\textrm{ref}i}-y_i(t)}{x_{\textrm{ref}i}-x_i(t)}\).
$$
To proceed, for any $k$ consider $T_k$ close to $T$ and the grid as in \eqref{grid} with $x_i^k:=x^k(t_i)$ for $i=1,\ldots,k$. Denote $t^k_0:=0$, $t^k_{i+1}:=\disp\sum_{j=0}^{i}h^k_j$, $I^k_i:=[t^k_i,t^k_{i+1})$ for $i=0,\ldots,k-1$, and $I^k_k:=\{T\}$. According to \eqref{Pi}, we have the algorithm
\begin{eqnarray}\label{vk}
\begin{array}{ll}
x^k_0\in A\;\mbox{ and }\;x^k_{i+1}:=x^k_i+h^k_iV^k_{i+1}\;\mbox{ with}\\
V^k_{i+1}:=\disp\Pi\big(g(x^k_i,u^k_i);C_{h^k_i}(x^k_i)\big),\;i=0,\ldots,k-1.
\end{array}
\end{eqnarray}
Taking into account the construction of $x^k_i$ for $0\le i\le k-1$, define a sequence of piecewise linear arcs $x^k\colon[0,T_k]\to\R^{2n}$, which pass through those points as
\begin{eqnarray}\label{xk}
x^k(t):=x^k_i+(t-t^k_i)V^k_{i+1}\;\mbox{ for all }\;t\in I^k_i,\;k\in\N.
\end{eqnarray}
We clearly have the relationships
\begin{eqnarray}\label{2m}
x^k(t^k_i)=x^k_i=\underset{t\to t^k_i}{\lim}x^k_i(t)\;\mbox{ and }\;\dot{x}^k_i(t):=V^k_{i+1}
\end{eqnarray}
for all $t\in(t^k_i,t^k_{i+1})$. As discussed in \cite{hb1}, the solutions to \eqref{xk} in the {\em uncontrolled} setting of \eqref{vk} with $g=g(x(t))$ uniformly converge on $[0,T_k]$ to a trajectory of a certain perturbed sweeping process. The {\em controlled} model under consideration here is significantly more involved. For all $x(t)\in\R^{2n}$, define the set
\begin{eqnarray}\label{K}
&K(x(t)):=\big\{y(t)\in\R^{2n}\big|\;D_{ij}(x(t))+\nonumber\\
&\nabla D_{ij}(x(t))(y(t)-x(t))\ge 0\;\mbox{ whenever }\;i<j\big\},
\end{eqnarray}
which allows us to represent the algorithm in \eqref{vk}, \eqref{xk} as
\begin{eqnarray*}
x^k_{i+1}=\Pi\big(x^k_i+h^k_ig(x^k_i,u^k_i);K(x^k_i)\big)\;\mbox{ for }\;i=0,\ldots,k-1.
\end{eqnarray*}
It can be equivalently rewritten in the form
\begin{eqnarray*}
&x^k\big(\vartheta^k(t)\big)=\Pi\big(x^k(\tau^k(t))+h^k_i g(x^k(\tau^k(t)),u^k(\tau^k(t));\\
&K(x^k(\tau^k(t))\big)
\end{eqnarray*}
for all $t\in[0,T_k]$ with $\tau^k(t):=t^k_i$ and $\vartheta^k(t):=t^k_{i+1}$ for all $t\in I^k_i$. Taking into account the construction of $K(x)$ in \eqref{K} together with \eqref{2m}, we arrive at the sweeping inclusions
\begin{eqnarray}\label{h:3.2}
&\dot x^{k}(t)\in-N\big(x^{k}(\vt^{k}(t));K(x^{k}(\tau^{k}(t)))\big)\nonumber\\
&+g\big(x^{k}(\tau^{k}(t)),u^{k}(\tau^{k}(t))\big)\;\mbox{ a.e. }\;t\in[0,T],
\end{eqnarray}
where $x^{k}(0)=x_0\in K(x_0)=A$ and $x^{k}(\vt^{k}(t))\in K(x^{k}(\tau^{k}(t)))$ on $[0,T]$. To formalize \eqref{h:3.2} as a controlled perturbed sweeping process, define the convex polyhedron
\begin{eqnarray}\label{e:131***}
C:=\bigcap\big\{x\in\R^{2n}\big|\;\la x^j_*,x\ra\le c_j,\;j=1,\ldots,n-1\big\}
\end{eqnarray}
with $c_j:=-(R_{USV}+R_{obs})$, where $R_{USV}$, $R_{obs}$ are the radii of the considered USV and the obstacle, respectively, and the $n-1$ vertices of the polyhedron
\begin{eqnarray*}\label{e}
x^j_*:=e_{j1}+e_{j2}-e_{(j+1)1}-e_{(j+1)2},\;j=1,\ldots,n-1,
\end{eqnarray*}
where $e_{kl}:=\big(\oe_{11},\oe_{12},\oe_{21},\oe_{22},\ldots,\oe_{k1},\oe_{k2},\ldots,\oe_{n1},\oe_{n2}\big)\in\R^{2n},$ $k=1,\ldots,n$ and $l=1,2$, with 1 at only one position of $\oe_{kl}$ and $0$ at all the other positions.

Let us formulate the {\em sweeping optimal control problem} (P) that can be treated as a continuous-time counterpart of the discrete algorithm to optimize the {\em controlled marine surface vehicle model}. Consider the cost functional
\begin{equation}\label{t:102*}
\mbox{minimize }\;J[x,u,T]:=\disp\frac{1}{2}\big\|x(T)\big\|^2,
\end{equation}
which reflects the model goal to minimize the distance and the time of the USV from the admissible configuration set to the target. We describe the continuous-time dynamics by the controlled sweeping process
\begin{align}\label{t:101*}
\left\{\begin{array}{lcl}-\dot{x}(t)\in N\big(x(t);C\big)+g\big(x(t),u(t)\big),\\
x(0)=x_0\in C,\;u(t)\in U\;\mbox{ a.e. }\;t\in[0,T],
\end{array}\right.
\end{align}
where $C$ is taken from \eqref{e:131***}, the control constraints reduce to \eqref{t:99*}, and the dynamic nonoverlapping condition $\|x^i(t)-x^j(t)\|\ge R_i+R_j$ is equivalent to the pointwise state constraints
\begin{eqnarray*}
x(t)\in C\Longleftrightarrow\la x^j_*,x(t)\ra\le c_j,\;t\in[0,T],\;j=1,\ldots,n-1.
\end{eqnarray*}

Now, we present the applications of the optimality conditions from Theorem~\ref{theorem:main thm} to the sweeping optimal control problem in \eqref{t:102*} and \eqref{t:101*} with two moving marine crafts MC 1 and MC 2. The marine surface vehicles are represented by triangle shapes immersed in discs (see Fig.~1). The objective is to move MC 1 and MC 2 to the target without colliding with each other. However, in the presence of MC 2, after the contacting time $t^*$ the vehicle MC 1 pushes MC 2 to the target with the same velocity. The mathematical USV's model is taken from the physical ship called Cyber-Ship \cite{s05} with the mass 23.8 kg and the length 1.255 m.

The initial configuration (positions of MC 1 and MC 2) is $x(0)=\(x^1(0),x^2(0)\)$, and the target is the origin. The radii of the discs used in this model are $R_1=R_2=3.5\;\textrm{m}$. Then, we have the model in \eqref{t:102*} and \eqref{t:101*} with the data
\begin{equation}\label{data}
\left\{
\begin{array}{ll}
n=2,\;x_*=\(1,1,-1,-1\),\;c=-7,\\[1ex]
g(x,u):=u,\;\vph(x,T):=\dfrac{1}{2}\|x(T)\|^2,\\[1ex]
U:=\big\{u=(u^1,u^2)\in\R^2\big|\;u^1\in [-2,2];\;u^2\in[-2,2]\big\},\\[1ex]
x^1(0)=\(-25,-25\),\;x^2(0)=\(-15,-15\),\\[1ex]
\th^u:=\th^u_1=\th^u_2=45^{\circ},\;s_1=s_2=1.
\end{array}\right.
\end{equation}
The set $C$ in \eqref{t:101*} is described by
\begin{align*}
C&=\big\{x\in\R^{4}:\;\la x_*,x\ra\le c\big\}\nonumber\\
&=\big\{x\in\R^{4}:\;x^{11}+x^{12}-x^{21}-x^{22}\le -7\big\}\nonumber\\
&=\big\{x\in\R^{4}:\;|x^{21}-x^{11}|+|x^{22}-x^{12}|\ge 7\big\}\nonumber\\
&\mbox{(under the imposed assumptions}\; x^{21}>x^{11}\mbox{ and }x^{22}>x^{12})\\
&=\big\{x\in\R^{4}:\;\|x^2-x^1\|\ge 2R\big\} \;\mbox{ for all }\;t\in[0,T].
\end{align*}
The structure of the problem suggests that the object only changes its velocity when it hits the boundary at some time $t_c$ with $c\in\{0,1,\ldots,k\}$. Moreover, if $t_c<T$, the object slides on the boundary of $C$ for the whole interval $[t_c,T]$. In this case, by construction of $t_c$ it must be one of the mesh points $t^k_i$ of some partition $\Delta_k$ in Theorem~\ref{theorem:main thm}. It is easy to see that all the assumptions of Theorem~\ref{theorem:main thm} are satisfied for \eqref{data}, and we can employ the obtained necessary optimality conditions, where the superscript ``$k$" is dropped, and where $\varrho_k,\;\(\xi_{iu},\;\xi_{iy}\)$ are supposed to be $0$ for large $k$ due to the convergence of discrete optimal solutions.
\begin{figure}[ht]
\centering
\includegraphics[scale=0.18]{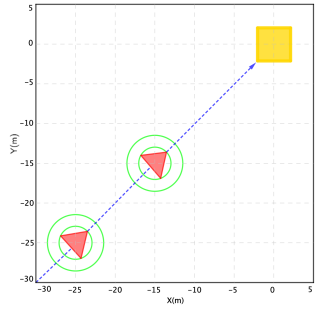}
\caption{Marine vehicle navigation before the contacting time.}
\label{Fig1}
\end{figure}
\begin{figure}[ht]
\centering
\includegraphics[scale=0.18]{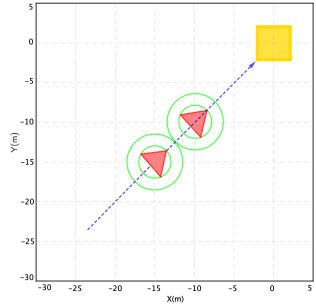}
\caption{Marine vehicle navigation after the contacting time.}
\label{Fig1}
\end{figure}
Applying all the conditions in Theorem~\ref{theorem:main thm} and using calculations in MATLAB, we get that the USV reaches the target at the minimum ending time $\overline T\approx 26.003$, and the hitting time is $t_c\approx 7.8201$. In this way we arrive at the optimal velocity $(\ou^1,\ou^2)\approx(1.67547,0.49999)$ and the optimal trajectory on $[0,\overline T]$ with the different expressions before and after the contacting time:
\begin{eqnarray*}
\left\{\begin{array}{ll}
\ox^1(t)\approx\(-25+1.18474t,-25+1.18474t\),\\
\ox^2(t)\approx\(-15+0.35355t,-15+0.35355t\),
\end{array}\right.
\end{eqnarray*}
for $t\in[0,7.8201)$, and
\begin{eqnarray*}
\left\{\begin{array}{ll}
\ox^1(t)\approx \(-21.7500+0.76914t,-21.7500+0.76914t\),\\
\ox^2(t)\approx \( -18.2500+0.76914t,-18.2500+0.76914t\),
\end{array}\right.
\end{eqnarray*}
for $t\in[7.8201,26.003]$.
\section{Concluding Remarks}\label{sec:Conclusions}

In this paper we formulated and studied a new class of optimal control problems governed by free-time controlled sweeping processes, where the duration of the process is also included to optimization. 
Developing the method of discrete approximations and using the generalized differential tools of variational analysis, we derive efficient necessary conditions for discrete optimal solutions that approximate a prescribed local minimizer of the continuous-time problem. The obtained results are applied to optimizing a controlled version of the marine surface vehicle model with static and dynamic obstacles, which is formulated in this paper based on the sweeping dynamics.

In our future research, we intend to furnish the limiting procedure of deriving necessary optimality conditions for the free-time continuous sweeping dynamics and provide further applications (qualitative and algorithmic) to more general versions of the controlled marine surface vehicle model dealing with many vessels.\vspace*{-0.05in}

\section*{Acknowledgment}
T.H. Cao acknowledges the support of the National Research Foundation of Korea grant funded by the Korea Government (MIST) NRF-2020R1F1A1A01071015.

N.T. Khalil and F.L. Pereira acknowledge the support of FCT R\&D Unit SYSTEC-POCI-01-0145-FEDER-006933 funded by ERDF—COMPETE2020—FCT/MEC—PT2020, Project STRIDE-NORTE-01-0145-FEDER-000033 funded by ERDF—NORTE 2020, and Project MAGIC-POCI-01-0145-FEDER-032485  funded by FEDER-COMPETE2020-POCI and PIDDAC through FCT/MCTES.

B.S. Mordukhovich, D. Nguyen, and T. Nguyen  acknowledges the support of the US National Science Foundation under grants DMS-1007132 and DMS-1512846, by the US Air Force Office of Scientific Research grant \#15RT0462.\vspace*{-0.05in}

\bibliographystyle{apalike}

\end{document}